# On the Extraction of Amicable Pairs Between Ibn Sīnā-al-Baghdādī and al-Kāshī

Mahmoud Annaby[1]

[1] Department of Mathematics, Faculty of Science, Cairo University,
Giza 12613, Egypt
mhannaby@@sci.cu.edu.eg

**Abstract.** In this note, we briefly analyse the works Ibn Sīnā (c.980-1037), al-Baghdādī (c.980-1037) and al-Kāshī (d. c. 1429) on the extraction of amicable pairs. We compare these works in the view of the well-known Thābit ibn Qurrā's rule. Contrary to the belief that these authors are merely stating variations of the famous Thābit ibn Qurrā's (d. 901) rule in different manners, we show that their statements are different. We prove that the statements of both ibn Sīnā and al-Baghdādī led to an unsolved conjecture, while al-Kāshī's statement led to the wrong amicable pair $(2024, 2296)$. The conjecture implied by ibn Sīnā and al-Baghdādī statements requires testing big primes corresponding to the known 51 Mersenne primes, and we show that no counterexamples yet exist.

**Keywords:** Amicable pairs, Mersenne primes, Fermat primes.

## 1 Introduction

Let $\sigma(n)$ denote the sum of proper divisors (aliquot parts) function

$$\sigma(n) := \sum_{p|n, p<n} p, \qquad n \in \mathbb{N} = \{1, 2, \cdots\}. \tag{1}$$

Two positive integers $m, n \in \mathbb{N}$ are said to be amicable if $\sigma(m) = n, \sigma(n) = m$. An example of amicable pairs is the smallest ever known amicable pair $(220, 284)$. According to Heath [6, p. 75], Iamblichus attributed this pair to Pythagoras (500 B.C). The first known rule to extract amicable pairs is the rule derived by Thābit b. Qurrā, [3, 10] and [12, p. 278]. Throughout this note we use the following notations for $n \in \mathbb{N}$:

$$\left. \begin{array}{l} a_n := 3 \times 2^{n-1} - 1, \quad b_n := 3 \times 2^n - 1, \quad c_n := 9 \times 2^{2n-1} - 1, \\ r_n := 2^n \times a_n \times b_n, \quad s_n := 2^n \times c_n, \qquad m_n := 2^n - 1. \end{array} \right\} \tag{2}$$

Thābit b. Qurra proved the following theorem, after proving nine preliminary lemmas.

**Theorem 1.** Let $n > 1$. If $a_n, b_n, c_n$ are prime, then $r_n, s_n$ are amicable.

This rule does not only establish sufficient conditions to extract amicable pairs, but it also gives the pairs in closed forms. The first known Pythagorean pair $(220, 284)$ corresponds the case $n = 2$. The second pair, which corresponds, $n = 4$, i.e. $(17296, 18416)$, is attributed to Fermat, and the third pair corresponding to $n = 7$, i.e. $(9363584, 9437056)$, is attributed to Descartes, who both rediscovered Thābit b. Qurra's rule between the years 1636-1638. However, Rashed in [11] showed that the pair $(17296, 18416)$ was calculated by al-Fārisī, much earlier than Fermat, namely before 1319, [12, p. 284-285], or even as earlier as 1307, [12, p. 313], or 1300 [3, 16, p. viii]. There are also some works that attribute al-Fārisī-Fermat pair to Thābit b. Qurra himself, an argument that requires more clues, [4, 7] and [13, p. 337]. Unlike, to the case of the second pair, al-Yazdī (d. 1637) claimed the discovery of the third pair in his *Treatise on Arithmetic*, [12, 286] and therefore the third pair $(9363584, 9437056)$ is known around the year 1600, [16, p. viii]. No other pairs are known due to b. Qurra's rule. Table 1 below presents the calculations associated with the three known amicable pairs derived according to Thābit ibn Qurrā's rule, see also [3,16]. Lemma 1 below proves that if $a_n, b_n$ are prime, then $r_n, s_n$ are amicable if and only if $c_n$ is prime.

Thābit b. Qurra's rule plays a central role in Arabic mathematical texts on the theory of numbers, cf. [11, pp. 209-278]. In many Arabic texts authors believe that it is possible to relax the sufficient conditions in two manners. The first one is to reduce the number of primes. The second approach is to replace the condition of the primality of $c_n$, which is necessary and sufficient for $r_n, s_n$ to be amicable, provided that $a_n, b_n$ are prime, by other less computationally expensive prime numbers like assuming that $m_{n+1}$ is prime, i.e., a Mersenne prime. We concentrate on the works of three Arabic authors, b. Sīnā (c.980-1037) [15], al-Baghdādī (c.980-1037) [2], and al-Kāshī (c. 1380 – 1429) [1, 8].

We summarize the works of these three authors and point out to the conjectures implied by their statements as well as the pitfall in the derivations of amicable pairs, if any. It is widely believed that all three scientists are just repeating the famous rule of b. Qurra to determinate amicable pairs. For instance, in a recent translation of al-Kāshī's *Miftāḥ al-ḥisāb*, (A Key for Arithmetic) [1], the authors translated al-Kāshī's work on amicable pairs and mentioned in their commentary that his rule, the fiftieth rule is ibn Qurra's [1, p. 121]. They did not notice that al-Kāshī's reduced ibn Qurra's sufficient conditions to two instead of three, which is wrong. In [1], it is not also recognized that this major mistake led al-Kāshī to claim that $(2024, 2296)$ is another amicable pair, which is not true since

$$\sigma(2024) = 2296, \sigma(2296) = 2744 \neq 2024. \tag{3}$$

However, while computing the aliquot parts of the amicable pairs, al-Kāshī assumed the primality of $c_n$. This pairs are sociable since 2024 is a father of 2296. Due to al-Baghdādī [2, p. 231], for $m, n \in \mathbb{N}$ $m$ is called a father of $n$ if $\sigma(m) = n$. Also, [12, 282], it is mentioned that b. Sīnā as well as all-Baghdādī statements are also slight variations of b. Qurra's theorem. In fact, both replaced the condition that $c_n$ is prime by assuming that $m_{n+1}$ is prime in addition to the other two conditions, i.e., the primality of both $a_n, b_n$ see [12, p. 312]. This leads to no known mistakes, but to an open conjecture. Nevertheless, the statements of b. Sīnā and al-Baghdādī are formally different, but they lead to the same conjecture. Additionally, while al-Baghdādī is constructing a calculation of the first amicable pair he used a technique that may lead to an algorithm, which we prove that it works only once. However, it is not easy to derive counter examples since the technique is connected to Fermat primes and all known Fermat primes do not lead to any counterexamples.

## 2 Amicable Pairs in Ibn Sīnā's "*al-Arithmāṭīqī*"

In [15], b. Sīnā derived his own rule by replacing that $c_n$ is prime, by the condition that $m_{n+1}$ is prime. In other words[1]:

> If you sum even-even numbers together with one and the sum is prime, for which if we add to them the last term and subtract its predecessor, the results after increment and decrement are prime, then multiplying the number with increment by the number with decrement, and multiplying the result by the last term, results a number with a lover; and his lover is the number which results from adding the number with a lover to the sum of the mentioned abundant and deficient numbers multiplied by the last term, and they are amicable.

The sum of even-even numbers together with one is

$$1 + 2 + 2^2 + \cdots + 2^n = 2^{n+1} - 1 = m_{n+1}. \tag{4}$$

Since:
$$\left.\begin{array}{l} m_{n+1} + 2^n = 2^{n+1} - 1 + 2^n = 3 \cdot 2^n - 1 = b_n, \\ m_{n+1} - 2^n = 2^{n+1} - 1 - 2^{n-1} = 3 \cdot 2^{n-1} - 1 = a_n, \end{array}\right\} \tag{5}$$

and that

---

[1] ibn Sina, Al-Shifāʾ, al-Ḥisāb (al-Arithmāṭīqī), [15, p. 28].

إذا جمعت أعداد زوج الزوج والواحد معهما فاجتمع عدد أول بشرط أن يكون إذا زيد عليهما آخرها ونقص الذي قبله كان المبلغ بعد الزيادة والمبلغ بعد النقصان أولياً فضرب المبلغ المزيد عليه في المبلغ المنقوص ثم ضرب ما اجتمع في آخر المجموعات حصل عدد له حبيب، وحبيبه العدد الذي يكون من زيادة مجموع الزائد والناقص المذكورين ضرباً في آخر المجموعات على العدد الموجود أولاً الذي له حبيب وهما متحابان.

$$\left.\begin{aligned} 2^n \times a_n \times b_n &= r_n, \\ 2^n(a_n + b_n) + r_n &= s_n, \end{aligned}\right\} \quad (6)$$

then b. Sīnā's conjecture can be stated as, [12, p. 312]:

**Conjecture 1.** Let $n > 1$. If $a_n, b_n, m_{n+1}$ are prime, then $r_n, s_n$ are amicable.

In the case $n = 7$, which led to al-Yazdī-Descartes pair $(9363584, 9437056)$ according to Theorem, indicates that while $a_n, b_n, c_n$ are primes, $m_{7+1} = 255$ is not prime. Therefore, b. Sīnā's conjecture is not implied by Theorem 1 above. The other direction is not as easy task as will be seen below. But at least, we see also from the case $n = 6$, that the primality of $m_{n+1}$ does not imply the primality of $c_n$ since $m_{6+1} = 2^7 - 1$ is prime, but $c_6 = 9 \times 2^{11} - 1$ is not prime.

The question now is to prove or disprove b. Sīnā's conjecture. We start by the following lemma.

**Lemma 1.** Let $a_n, b_n, m_n$ be prime. If $c_n$ is not prime, then $r_n, s_n$ are not amicable.

**Proof.** Since $c_n$ is composite, then direct computations yield

$$\sigma(s_n) = \sigma(2^n c_n) > \sum_{k=0}^{n-1} 2^k c_n + \sum_{k=0}^{n} 2^k = 2^n a_n b_n = r_n.$$

∎

Lemma 1 does not prove that Conjecture 1 is wrong. Consequently, Conjecture 1 is wrong if we can find $n \in \mathbb{N}, n > 1$, for which $a_n, b_n, m_{n+1}$ are prime, but $c_n$ is not prime, or, in other words, to prove the following statement:

$$\exists n > 1 \ (a_n \text{ is prime} \land b_n \text{ is prime} \land m_{n+1} \text{ is prime} \land c_n \text{ is not prime}). \quad (7)$$

Thus, a disprove of Conjecture 1 requires a counterexample. So, we first seek a Mersenne prime $m_{n+1}$ and then we check if $a_n, b_n$ are also prime, i.e. the sufficient conditions of Conjecture 1 are fulfilled. If $a_n, b_n$ are also prime, then a disprove requires us to prove that the corresponding $c_n = a_n + b_n + a_n b_n$ is not prime. Since there are only known 51 Mersenne primes[2], we can check them all to seek a desired counterexample. Using the Computer Algebra System *Mathematica*, we found that no known counterexamples disprove Conjecture 1 as Table 3 below proves. Although $a_n, b_n$ as well as $m_{n+1}$ are big numbers the primality is checked via simple theorems like the fact stated in [9, p. 445], which states that $3 \times 2^n - 1$ is not prime if $n \equiv 1 \pmod 4$, since in this case 5 will be a factor of $3 \times 2^n - 1$.

---

[2] https://www.mersenne.org/primes/

## 3   Amicable Pairs in al-Baghdādī's "*al-Takmila*"

Despite the fact proved below, that al-Baghdādī statement leads to b. Sīnā's conjecture to construct amicable pairs, his approach is different. He started by deriving the first pair $(220, 284)$ and then to give the general rule. However, the way he derived the first pair is interesting. In his statement[3]:

> The genesis in collecting amicable numbers is from even-even numbers: So, we start with the first even, which is four, we add one to it, it becomes five, which is a prime number, then we multiply this five by two, and we add one to the result, it becomes eleven, and as the eleven is prime, and if we subtract four from it, the remainder will be prime: it is multiplied by five, then by four, which is the first even-even, results two hundred and twenty, which is one of the two amicable numbers. And if we want the second number, we gather between the five and the eleven, it becomes sixteen, which we multiply by the four, which is the first even-even, it becomes sixty four, then we add it to the first number, which is two hundred and twenty, it becomes two hundred and eighty four, and it is the second number.

The above description is made by al-Baghdādī only to establish the first amicable pair because al-Baghdādī is using the prime numbers 5, 11 as initial values in his calculations of the general rule. Now let us state his example in a general manner to extract amicable pairs. It is worthwhile to mention that he did not state the general setting we stated below. We will mention his general statement, that coincide with b. Sīnā's conjecture later. Now let us demonstrate the general view he used to derive the first example. For $n \geq 2$, let $\alpha_n, \beta_n, \gamma_n$ be prime, where

$$\begin{aligned}
\alpha_n &= 2^n + 1, \\
\beta_n &= 2(2^n + 1) + 1 = 2^{n+1} + 3, \\
\gamma_n &= \beta_n - 2^n = 2^n + 3.
\end{aligned}$$

Then the following numbers are amicable

$$\begin{aligned}
\lambda_n &= 2^n \alpha_n \beta_n = 2^{3n+1} + 5 \cdot 2^{2n} + 3 \cdot 2^n, \\
\mu_n &= 2^n(\alpha_n + \beta_n) + 2^n \alpha_n \beta_n = 2^{3n+1} + 8 \cdot 2^{2n} + 7 \cdot 2^n.
\end{aligned}$$

---

[3] al-Baghdādī, *al-Takmila*, [2, p. 230].

والأصل في تحصيل الأعداد المتحابة، من أعداد زوج الزوج: فنبدأ بأول زوج وهو أربعة، ونزيد عليه واحداً، فيصير خمسة، وهي عدد أولي، فنضرب هذه الخمسة في اثنين، ونزيد على المبلغ واحداً، فيصير أحد عشر، فإن كان الأحد عشر أولياً، وإذا نقصنا منه الأربعة كان الباقي أولياً: ضرب في الخمسة، ثم في الأربعة التي هي أول زوج زوج، فبلغ مئتين وعشرين، وهو أول العددين المتحابين. فإن أردنا العدد الثاني، جمعنا بين الخمسة والأحد عشر فيصير ستة عشر، فضربناها في الأربعة، وهي أول زوج زوج، فيصير أربعة وستين، فنزيدها على العدد الأول، وهو مئتان وعشرون، فيبلغ مئتين وأربعة وثمانين، وهي العدد الثاني.

If $n = 2$, then
$$\alpha_2 = 5, \quad \beta_2 = 11, \quad \gamma_2 = 7, \quad \lambda_2 = 220, \quad \mu_2 = 284.$$

The first conjecture implied by al-Baghdādī's derivation of the first amicable pair is the following.

**Conjecture 2.** Let $n > 1$. If $\alpha_n, \beta_n, \gamma_n$ are prime, then $\lambda_n, \mu_n$ are amicable.

There is no connection between Conjecture 2 and Conjecture 1. Notice for example that, $\lambda_n, \mu_n$ are different from $r_n, s_n$, which are:
$$\begin{aligned} r_n &= 9 \cdot 2^{3n-1} - 9 \cdot 2^{2n-1} + 2^n, \\ s_n &= 9 \cdot 2^{3n-1} - 2^n. \end{aligned}$$

There are some clues that make us believe that he has conjecture 2 in mind while searching for the first amicable pair. Some of these are:

- Unlike b. Sīnā, who directly gave $(220, 284)$ as an example, al-Baghdādī stated his search in the form of an algorithm.
- While searching for $(220, 284)$ the primality of $\gamma_2$, i.e. 7, is required, although this condition has no role in deriving the first amicable pair, since 7 is neither a proper divisor for 220, nor for 284.
- Another reason is that he required the primality of a Fermat number $\alpha_2$, i.e. 5, while in his general algorithm Fermat primes did not play any role, but Mersenne primes have a role as in the case of b. Sīnā's conjecture as we will see below. Since the only known Fermat primes $f_n = 2^n + 1$, are those numbers, which correspond [14],

$$n = 0, 1, 2, 4, 8, 16,$$

and checking these cases we see, as Table 4 indicates that $\alpha_n, \beta_n, \gamma_n$ are prime only if $n = 2$, it is not easy to derive a counterexample to indicate that Conjecture 2 is wrong. However, direct computations indicate that setting $\sigma(\lambda_n) = \mu_n$, leads to the equation
$$2^{3n+1} + 9 \cdot 2^{2n} + 5 \cdot 2^n - 8 = 2^{3n+1} + 8 \cdot 2^{2n} + 7 \cdot 2^n,$$
which leads to the equation
$$2^{2n} - 2^{n+1} - 8 = 0,$$
and the latter occurs when and only when $n = 2$.

It remains probable that al-Baghdādī had not this conjecture in mind while deriving the first amicable pair, whose calculations are used as initial step in his iteration algorithm to derive other amicable pairs. Would he use Conjecture 2 to construct the first amicable pair, if the first computations did not lead to $(220, 284)$, or rather he was going to implement his general algorithm described above, which coincide with b. Sīnā's conjecture as we will see?

Let us now state the general algorithm of al-Baghdādī and prove that it is merely b. Sīnā's conjecture. He states that[4]:

> If we want what is after that from the amicable numbers, we multiply the eleven by two, and the five by two, and we add one to each of the results, and we check: if each number of them is prime, and if we subtract the even-even number that succeeds what we have subtracted first, the remainder becomes a prime number, we multiply one by the other, and then the result by the second even-even, which is eight. And if they are not from prime numbers, there will be no two amicable numbers from them, and they have to be multiplied by two and one is added to each one of the results, and we continue doing this until the matter ends with them to be from the prime numbers, then we multiply one by the other, and then the result by the even-even in its rank, and the result is one of the two amicable numbers. And this is its measure.

In mathematical words, al-Baghdādī says that: If
$$A_n = 5 \cdot 2^n + 2^{n-1} + 2^{n-2} + \cdots + 1,$$
$$B_n = 11 \cdot 2^n + 2^{n-1} + 2^{n-2} + \cdots + 1,$$
$$M_n = B_n - 2^{n+2}$$
are prime, then
$$R_n := 2^n \times A_n \times B_n, \qquad S_n := 2^n(A_n + B_n) + R_n$$
are amicable. If we simplify this, we find that:
$$A_n = 3 \cdot 2^{n+1} - 1 = a_{n+2},$$
$$B_n = 3 \cdot 2^{n+2} - 1 = b_{n+2},$$
$$M_n = 2^{n+3} = m_{n+3}.$$

Thus $R_n = r_{n+2}$, $S_n = s_{n+2}$, and al-Baghdādī's statement coincide entirely with b. Sīnā's one with a shift, and there is no need to repeat the arguments mentioned above.

It remains an open question to ask why both of them replaced the primality of $c_n$ by the primality of $m_{n+1}$? Did they have another common source? Are both authors communicated? There are some arguments that these ideas had not been exchanged. First, both statements are different, although they lead to the same conjecture. Sec-

---

[4] al-Baghdādī, *al-Takmila*, [2, pp. 230-231]:

فإن أردنا ما وراء ذلك من الأعداد المتحابة، ضربنا الأحد عشر في اثنين، والخمسة في اثنين، وزدنا على كل واحد من المبلغين واحداً، ونظرنا: فإن كان كل واحد مهما أولياً، وإذا أسقط زوج الزوج الذي يلي ما أسقط أولاً، بقي عدد أولي، ضرب أحدهما في الآخر، ثم ما بلغ في زوج الزوج الثاني، وهو ثمانية. وإن لم يكونا من الأعداد الأولية، فلا يكون منهما عددان متحابان، وينبغي أن يضربا في اثنين ويزاد على كل واحد من المبلغين واحد، ثم لا يني يفعل ذلك إلى أن ينتهي الأمر فيهما إلى أن يكون العددان من الأعداد الأولية، فنضرب أحدهما في الآخر، ثم ما بلغ في زوج الزوج الذي في رتبته، فما بلغ فهو أحد العددين المتحابين. وهذا قياسه.

ondly, al-Baghdādī stated in [2] clearly that there is no perfect numbers between $10^4$ and $10^5$[5], but b. Sīnā, [15] stated that there is a perfect number in each interval $10^{n-1} < x \leq 10^n, n \in \mathbb{N}$[6].

## 4  Amicable Pairs in al-Kāshī's "*Miftāḥ al-Ḥisāb*"

In his book "*Miftah al-Hisāb*" (A Key for Arithmetic) al-Kāshī introduced his own rule to extract amicable pairs[7] as follows:

> If we want to extract the two amicable numbers, which are two numbers such that the sum of the aliquot parts of each of them is equal to the other, then we look for a number from the powers of two such that if we multiply it once by one and a half and by three another time, and subtract one from both results in which no number other than the one counts the remainders, if it exists, we call the first remainder the first odd and the second the second odd. The second odd necessarily exceeds the double of the first odd by one, then, we multiply the first odd by the second odd, and we call the result the third odd, then, we multiply the existing number from the powers of two once by the third odd and another by the sum of the first and the second odds, so the first is one of the two amicable numbers, and if we add the second result to it, then whatever it is, it will be the last one of the two amicable numbers.

Comparing the above translation with [4, p. 121], we have made more appropriate terms. For instance, we used word *power* instead of the word *multiple*, as it is more appropriate with the sequel. In fact, al-Kāshī used the word (تضاعيف) for powers, as

---

[5] al-Baghdādī, *al-Takmila*, [2, p. 227].

[6] ibn Sina, Al-Shifāʾ, al-Ḥisāb (al-Arithmāṭīqī), [15, p. 32]. In fact his statement is confusing because he says: "Be informed that the perfect number is necessarily even because it is constructed from the product of an odd number by an even, and it is agreed that what lies from it in the digits one which is six, and in the tenth one which is twenty eight, and in the hundreds one which is four hundred ninety six, and in the thousands one which is eight thousand one hundred and twenty, and continually one in every kind, with digits that are six or eight, although experimentally the alternation between them is not necessary" Therefore if he knows that the 6 and 8 do not necessarily alternate, he should have known the fifth and sixth perfect numbers which both end with 6, namely 33550336 and 8589869056. From the same fact he should have known that there are no perfect numbers between $10^4$ and $10^5$.

[7] In *Miftāḥ al-Ḥisāb*, p. 223, al-Kāshī stated his rule as follows:

إذا أردنا أن نستخرج العددين المتحابين وهما عددان يكون مجموع أجزاء كل واحد منهما مساوياً للآخر ، نطلب عدداً من تضاعيف الاثنين، إذا ضربناه تارة في واحد ونصف، وتارة في ثلاثة، وننقص من كل واحد من الحاصلين واحداً، فلا يعد لكل واحد من الباقيين غير الواحد، فإذا وجد يسمى الباقي الأول الفرد الأول، والثاني الفرد الثاني. ولابد أن يكون الفرد الثاني زايداً على ضعف الفرد الأول بواحد، ثم نضرب الفرد الأول في الفرد الثاني، ونسمي الحاصل بالفرد الثالث، ثم نضرب العدد الموجود من تضاعيف الاثنين تارة في الفرد الثالث وتارة في مجموع الفردين الأول والثاني، فيكون الأول أحد العددين المتحابين و إذا نزيد الحاصل الثاني عليه فما بلغ فهو العدد الأخير من المتحابين.

indicated in [8, p. 209][8]. Let us now translate Thus al-Kāshī's rule into mathematical language. We start with a number $2^n$, then we define
$$\frac{3}{2} \times 2^n - 1, \quad 3 \times 2^n - 1,$$
which are merely $a_n$ and $b_n$ respectively. And if $a_n$ and $b_n$ are prime, then
$$2^n \times (a_n \, b_n), \quad 2^n \times (a_n + b_n) + 2^n \times (a_n \, b_n),$$
are amicable. Noting that these last two numbers are merely $r_n$ and $s_n$, al-Kāshī's conjecture can be written as follows.

**Conjecture 3.** If $a_n, b_n$ are prime, then $r_n, s_n$ are amicable.

It is readily seen that this conjecture is wrong since from Table 1, we see that when $n = 3$, both $a_3 = 11, b_3 = 23$ are prime, but $r_3 = 2024, s_3 = 2296$ are not amicable as we indicated above. al-Kāshī did not assume primality of $c_n$, nor he assumed primality of the corresponding Mersenne prime as in b. Sīnā's- al-Baghdādī's conjecture. We may also believe that he has forgotten both assumptions, or it is missed from manuscripts, but this is not the case as he drived the wrong amicable pair $(2024, 2296)$ in his book. We have seen from Equation (8) above that $\sigma(2024) = 2296$, i.e. 2024 is a father of 2296. As a matter of fact, this always happens if $a_n, b_n$ are prime according to the following lemma, no matter $c_n$ is prime or not.

**Lemma 2.** Let $a_n, b_n$ be prime. Then $r_n$ is father of $s_n$, i.e. $\sigma(r_n) = s_n$.

**Proof.** Since $a_n, b_n$ be prime, it follows that

$$\begin{aligned} \sigma(2^n a_n b_n) &= \sum_{k=0}^{n} 2^k + \sum_{k=0}^{n} 2^k a_n + \sum_{k=0}^{n} 2^k b_n + \sum_{k=0}^{n-1} 2^k a_n b_n \\ &= \sum_{k=0}^{n} 2^k (1 + a_n + b_n) + \sum_{k=0}^{n-1} 2^k a_n b_n. \end{aligned} \qquad (9)$$

But

$$\begin{aligned} 1 + a_n + b_n &= 1 + 3 \cdot 2^{n-1} - 1 + 3 \cdot 2^n - 1 \\ &= 9 \cdot 2^{n-1} - 1, \end{aligned} \qquad (10)$$

and

$$\begin{aligned} a_n b_n &= (3 \cdot 2^{n-1} - 1)(3 \cdot 2^n - 1) \\ &= 9 \cdot 2^{2n-1} - 9 \cdot 2^{n-1} + 1. \end{aligned} \qquad (11)$$

Combining Equations (10) and (11) with Equation (9) yields

$$\sigma(2^n a_n b_n) = (2^{n+1} - 1)(9 \cdot 2^{n-1} - 1) \qquad (12)$$

---

[8] Ibid, p. 223, al-Kāshī defined in the ninth rule sums of the powers of the one, which are numbers that can be halved iteratively until, one reaches 2. See also [4, p. 79].

$$+(2^n - 1)(9 \cdot 2^{2n-1} - 9 \cdot 2^{n-1} + 1)$$
$$= 9 \cdot 2^{2n} - 9 \cdot 2^{n-1} - 2^{n+1} + 1$$
$$9 \cdot 2^{3n-1} - 9 \cdot 2^{2n-1} - 9 \cdot 2^{2n-1} + 9 \cdot 2^{n-1} + 2^n - 1.$$

Since $9 \cdot 2^{2n-1} + 9 \cdot 2^{2n-1} = 9 \cdot 2^{2n}$, and $2^{n+1} - 2^n = 2^n$, then (12) drives us directly to

$$\sigma(2^n a_n b_n) = 9 \cdot 2^{3n-1} - 2^n$$
$$= 2^n (9 \cdot 2^{2n-1} - 1) = s_n. \qquad (13)$$

This completes the proof. ∎

Throughout his computations, al-Kāshī derived the first pair $(220, 284)$ and the wrong pair $(2024, 2296)$, [8, pp. 224-225] and [1, pp. 121-125]. In addition, in [8, pp. 223-225] and [1, pp. 123-125] he stated all aliquot parts of both amicable pairs, which for $r_n = 2^n \times (a_n b_n)$ are:

$$1, 2, \cdots, 2^n, a_n, 2a_n, \cdots, 2^n a_n, b_n, 2b_n, \cdots, 2^n b_n, a_n b_n, 2a_n b_n, \cdots, 2^{n-1} a_n b_n,$$

and for $s_n = 2^n \times (a_n + b_n + a_n b_n) = 2^n c_n$ the aliquot parts are

$$1, 2, \cdots, 2^n, 1, 2, \cdots, 2^n, c_n, 2c_n, \cdots, 2^{n-1} c_n.$$

Thus, it is understood that he assumed primality of $c_n$, but he did not state it explicitly. In addition, he derived the wrong pair as he did not notice that $c_3 = 287$ is not prime and the table for aliquot parts [8, p. 225] and [1, p. 125] is wrong. This mistake has not been recognized neither in the commentaries of [1, pp. 121-125] nor in [8, pp. 223-225, pp. 322].

## 5 Appendix

This appendix involves tables of computations connected to the theorems and conjectures discussed above. The first table contains the computations of the first three amicable pairs due to Thābit b. Qurrā's rule. Notice that the case $n = 3$ leads to no amicable pairs. But we add it because we will refer to it.

**Table 1.** Amicable pairs extracted via Thābit b. Qurrā's rule with their discoverers.

| $n$ | $a_n$ | $b_n$ | $c_n$ | $r_n$ | $s_n$ | Discoverer |
|---|---|---|---|---|---|---|
| 2 | 5 | 11 | 71 | 220 | 284 | Pythagoras (500 B.C.) |
| 3 | 11 | 23 | 287 | 2024 | 2296 | No Pairs exist |
| 4 | 23 | 47 | 1151 | 17296 | 18416 | Fermat (1636) al-Fārisī (d.1320) |
| 7 | 191 | 383 | 73727 | 9363584 | 9437056 | Descatres (1638) al-Yazdī (d.1637) |

The next table, Table 2, is devoted to computations connected to b. Sīnā's-al-Baghdādī's conjecture. As we see only two pairs are implied by this conjecture. Notice that amicable pairs arose in the cases $n = 2, 4$ and that the rule is not applicable for the third pair $(9363584, 9437056)$ since $m_8 = 255$ is not prime.

**Table 2.** Amicable pairs according to b. Sīnā's- al-Baghdādī's conjecture.

| $n$ | $a_n$ | $b_n$ | $c_n$ | $m_{n+1}$ | $r_n$ | $s_n$ |
|---|---|---|---|---|---|---|
| 2 | 5 | 11 | 71 | 7 | 220 | 284 |
| 3 | 11 | 23 | 287 | 15 | 2024 | 2296 |
| 4 | 23 | 47 | 1151 | 31 | 17296 | 18416 |
| 7 | 191 | 383 | 73727 | 255 | 9363584 | 9437056 |

Table three is devoted to search for a counterexample to disprove b. Sīnā-al-Baghdādī's conjecture. We checked the primality of the three numbers $a_n, b_n, m_{n+1}$ through all known Mersenne primes and there is no known counterexample. As we see, no known examples exit to fulfil statement (7). The case $n = 1$ is set only to complete the table.

**Table 3.** An examination of all 51 known Mersenne primes $m_n$ to find a counterexample that disproves Sīnā- al-Baghdādī's conjecture.

| $n$ | $m_{n+1}$ is Prime | $a_n \wedge b_n$ are Prime | $n$ | $m_{n+1}$ is Prime | $a_n \wedge b_n$ are Primes | $n$ | $m_{n+1}$ is Prime | $a_n \wedge b_n$ are Primes |
|---|---|---|---|---|---|---|---|---|
| 1 | T | T | 3216 | T | F | 1398268 | T | F |
| 2 | T | T | 4252 | T | F | 2976220 | T | F |
| 4 | T | T | 4422 | T | F | 3021376 | T | F |
| 6 | T | F | 9688 | T | F | 6972592 | T | F |
| 12 | T | F | 9940 | T | F | 13466916 | T | F |
| 16 | T | F | 11212 | T | F | 20996010 | T | F |
| 18 | T | F | 19936 | T | F | 24036582 | T | F |
| 30 | T | F | 21700 | T | F | 25964950 | T | F |
| 60 | T | F | 23208 | T | F | 30402456 | T | F |
| 88 | T | F | 44496 | T | F | 32582656 | T | F |
| 106 | T | F | 86242 | T | F | 37156666 | T | F |
| 126 | T | F | 110502 | T | F | 42643800 | T | F |
| 520 | T | F | 132048 | T | F | 43112608 | T | F |
| 606 | T | F | 216090 | T | F | 57885160 | T | F |
| 1278 | T | F | 756838 | T | F | 74207280 | T | F |
| 2202 | T | F | 859432 | T | F | 77232916 | T | F |
| 2280 | T | F | 1257786 | T | F | 82589932 | T | F |

Table 4 checks correctness of Conjecture 2 over the known Fermat primes and as we see no counterexamples exist. However, this proposed conjecture is proved above to be wrong. Nevertheless, the conjecture is not proposed by al-Baghdādī explicitely, but, due to the above-mentioned reasons it might be in his mind while searching for the amicable pair $(220, 284)$.

Table 4. Computations of the primality of $\alpha_n, \beta_n, \gamma_n.$

| $n$ | $\alpha_n$ is Prime | $\beta_n \wedge \gamma_n$ are Prime |
|---|---|---|
| 2 | T | T |
| 4 | T | F |
| 8 | T | F |
| 16 | T | F |